\documentstyle{amsppt}
\magnification=\magstep1
\hcorrection{0in}
\vcorrection{0in}
\pagewidth{6.3truein}
\pageheight{9truein}
\NoRunningHeads

\def\Q{\Bbb Q}

\def\C{\Bbb C}
\def\P{\Bbb P}
\def\O{\Cal O}
\def\rup#1{\left\lceil #1 \right\rceil}
\def\rdown#1{\left\lfloor #1 \right\rfloor}
\def\Bs#1{\mathop{\text{\rm Bs}}\left| #1 \right|}
\def\Supp{\mathop{\text{\rm Supp}}}

\def\mult{\mathop{\text{\rm mult}}\nolimits}

\def\exc{{\text{exc}}}
\def\dmin{\delta_{\text{min}}}

\font\maru=lcircle10
\topmatter
\title Effective base point freeness on normal surfaces
\endtitle
\author Takeshi Kawachi \endauthor
\address
Department of Mathematics, Tokyo Institute of Technology\hfil\break
2-12-1 Oh-okayama, Meguro-ku, Tokyo 152-8551, JAPAN\hfil\break
Telefax: +81-3-5734-2738 
\endaddress
\email kawachi\@math.titech.ac.jp \endemail
\abstract
{\obeylines
In almost all situations, some special things are happen on a singularity.
This specialty sometimes causes something unpleasant.
Generally, smooth is reviewed fine and singular nasty,
But in some area of geometry, it would be different.
For global generations, a base locus has malignant nature,
It often appears at a smooth point, at singular 'tis rare.
Thus we would say ``fair is foul and foul is fair'' as in the classic literature,
So we would ``hover through the fog and filthy air''
To proceed the theory with no harm
And in caution against the witches charm.}
\endabstract
\endtopmatter
\document
\subhead Contents\endsubhead
\par
0. Introduction.\par
1. The invariants for singularities.\par
2. The main theorem.\par
3. A proof of the main theorem.\par
\par
\subhead Notation\endsubhead
\halign{\strut #\hfil & #\hfil \cr
  $\rup{\cdot}$   & the round up \cr
  $\rdown{\cdot}$ & the round down \cr
  $\{\cdot\}$     & the fractional part \cr
  $f^{-1}D$       & the strict transform (proper transform) of $D$\cr
  $f^*D$          & the pull back (total transform) of $D$\cr
  $\equiv$        & numerical equivalence. \cr
  $\sim$          & linear equivalence.\cr
}
\par
\head 0. Introduction\endhead
Let $Y$ be a compact normal two-dimensioinal algebraic space over $\C$
(``normal surface'' for short).
Let $B$ be an effective $\Q$-Weil divisor on $Y$ such that $\rdown{B}=0$.
Let $y\in Y$ be a given point, and $M$ be a nef and big $\Q$-Weil divisor
such that $K_Y+M+B$ is Cartier.
Various numerical conditions on $M$ which gives $y\not\in\Bs{K_Y+M+B}$
were studied as [ELM], [KM] and [F].
The following theorem unifies all earlier results and gives the first
effective version if $y$ is a log-terminal singularity
on a log-surface $(Y,B)$.
\par
\proclaim{Theorem 1}
Let $Y$, $y\in Y$, $M$ be as in the last paragraph.
Let $\delta$ and $\delta'$ be the invariant defined below.
Assume that $M^2>\delta$ and $M\cdot C\geq\delta'$
for any irreducible curve $C$ on $Y$ passing through $y$.
(We use the Mumford's $\Q$-valued pull-back and intersection theory
for $\Q$-Weil divisors on normal surfaces).
Then we have $y\not\in\Bs{K_Y+\rup{M}}$
\endproclaim
\par
\subhead Definition of $\delta$ and $\delta'$ \endsubhead
\par
Let $B=\rup{M}-M$.
Let $f\:(X,f^{-1}B)\to (Y,B)$ be the minimal desingularization of the germ
$(Y,y)$ if $y$ is singular, be the blowing up of $Y$ at $y$ if $y$ is smooth
(See Remark (1) below).
Let $\Delta_B=f^*(K_Y+B)-(K_X+f^{-1}B)$ be the canonical cycle and
$Z$ be the fundamental cycle of $y$.
Note that $\Delta_B$ is effective $\Q$-divisor supported on $f^{-1}(y)$
if $y$ is singular; $\Delta_B=(\mult_y B-1)Z$ if $y$ is smooth.
Note also that $\Delta_B$ is the negative of the usual discrepancy.
\par
\definition{Definition}
$$
  \dmin=
    \{-(Z-\Delta_B+x)^2\mid\text{$x$ is an effective $\Q$-Weil divisor
      supported on $f^{-1}(y)$}\}
$$
$$
\delta=
    \cases
      \dmin, 
         &\text{if $(Y,B)$ is log-terminal at $y$ (See Remark (2) below).} \\
      0, &\text{if $(Y,B)$ is not log-terminal at $y$.} \\
    \endcases
$$
\enddefinition
\par
Let $\Delta_B=\sum e_iE_i$ be the prime decomposition.
Since all log-terminal singularities on surfaces is classified by
[Al] and [Ky], we give the following definition.
\par
\definition{Definition}
$$
  \delta'=
    \cases
      1-\max\{e_1,e_n\}, &
        \text{if $(Y,B)$ is log-terminal at $y$ of type $A_n$,} \\
      & \text{where $E_1$ and $E_n$ are placed on the edge} \\
      & \text{of the chain of the dual graph} \\
      \text{any positive number}, &
        \text{if $(Y,B)$ is log-terminal at $y$ of type $D_n$} \\
      0, & otherwise
    \endcases
$$
\enddefinition
\par
Theorem 1 is useful because one has the folloing bound on
$\delta_{B,y}=-(Z-\Delta_B)^2$.
\par
\proclaim{Proposition 2}
\roster
\item $\delta_{B,y}\leq 4$ if $y$ is smooth point;
\item $\delta_{B,y}\leq 2$ if $y$ is a rational double point (RDP for short);
\item $\delta_{B,y}<2$ if $y$ is log-terminal but not smooth or an RDP;
\endroster
\endproclaim
\par
\remark{Remarks}
(1) The resolution of singularities $f\:(X,f^{-1}B)\to(Y,B)$
is only the desingularization $f\:X\to Y$ of $Y$.
So $f^{-1}B$ may have singularities on $f^{-1}(y)$,
or $f^{-1}B$ may not be normal crossings with the exceptional locus
$f^{-1}(y)$.
\par
(2) We are using non-standard definition of log-terminal.
we say $y$ is a log-terminal singularity on $(Y,B)$
if $[B]=0$ and all coefficients in the expression of $\Delta_B$ are
strictly less than 1,
irrespective of where $K_Y+B$ is $\Q$-Cartier.
\par
(3) If $y$ is at worst an RDP then Theorem 1 is essentially
the theorem in [EL].
If $B=0$ then Theorem 1 is essentially the theorem in [KM].
If $y$ is not log-terminal, the result is proved in [ELM].
In the case that $y$ is log-terminal, the minimality of $\dmin$ may give
smaller numerical conditions.
\endremark
\par
In section 1, we recall several definitions and the properties
concerning of $\delta_{B,y}$.
\par
In section 2, we recall Theorem 1 and introduce its corollaries.
\par
In section 3, we prove Theorem 1
\par
\head 1. The invariants for singularities\endhead
\subhead 1.1\endsubhead
Let $Y$ be a complete normal algebraic surface.
Let $f:X\to Y$ be a resolution of singularities of $Y$.
We use Mumford's $\Q$-valued pullback and intersection theory on $Y$
(cf\. [Mu]).
\par
Let $D$ be a $\Q$-Weil divisor on $Y$.
We write $f^*D=f^{-1}D+D_\exc$,
where $f^{-1}D$ is the strict transform of $D$ and
$D_\exc$ is the $f$-exceptional $\Q$-divisor on $X$ such that
$f^*D\cdot E=0$ for every $f$-exceptional curve $E$.
Since the intersection matrix of the $f$-exceptional curves
is negative definite, $D_\exc$ is uniquely determined.
Also the negative definiteness gives $f^*D\geq 0$ if $D\geq 0$.
\par
\definition{Definition}
Let $M$ be a $\Q$-Weil divisor on $Y$.
Then $M$ is {\it nef\/} if $M\cdot C\geq 0$ for any irreducible curve $C$
on $Y$.
Assume $M$ is nef,
then $M$ is {\it big} if in addition $M^2>0$.
\enddefinition
\par
\subhead 1.2\endsubhead
Let $y$ be a fixed point on $Y$.
Let $f\:X\to (Y,y)$ be the minimal resolution of the germ $(Y,y)$
if $y$ is singular, be the blowing up of $Y$ at $y$ if $y$ is smooth.
Let $f^{-1}(y)=\cup_{j=1}^{n} E_j$
\par
\definition{Definition}
The {\it fundamental cycle} of $(Y,y)$ is the smallest nonzero
effective divisor $Z=\sum z_jE_j$ such that $Z\cdot E_j\leq 0$ for all $j$
(cf\. [Ar]).
\enddefinition
\par
Since $f^{-1}(y)$ is connected, $z_j\geq 1$ for all $j$.
If $y$ is smooth, then $Z=E_1$ is the exceptional $(-1)$-curve.
Let $p_a(Z)=\frac12 Z\cdot(K_X+Z)+1$.
Since [Ar, Theorem 3], $p_a(Z)\geq 0$.
\par
\subhead 1.3\endsubhead
Let $\Delta=f^*K_Y-K_X$ be the $f$-exceptional $\Q$-Weil divisor.
Since $\Delta$ is determined by the intersection numbers,
$\Delta$ is uniquely defined as in (1.1).
Note that $\Delta\cdot E= -K_X\cdot E$ for any $f$-exceptional divisor $E$
by its definition.
\par
If $y$ is smooth, then $\Delta=-E_1$.
If $y$ is singular, $\Delta$ is effective since $f$ is the minimal resolution.
\par
\definition{Definition}
$\Delta$ is called as the {\it canonical cycle} of $(Y,y)$.
\enddefinition
\par
\subhead 1.4\endsubhead
Let $B=\sum b_iC_i$ be an effective $\Q$-Weil divisor on $Y$
such that $b_i\leq 1$.
Let $\Delta_B=f^*(K_Y+B)-(K_X+f^{-1}B)$ be the $f$-exceptional
$\Q$-Weil divisor.
$\Delta_B$ is defined uniquely as in (1.3).
Let $\Delta_B=\sum e_jE_j$.
\definition{Definition}
$(Y,B)$ is {\it log-terminal} (resp. {\it log-canonical\/}) at $y$ if
\roster
\item $[B]=0$ (resp. $\rup{B}$ is a reduced divisor),
\item $e_j<1$ (resp. $e_j\leq 1$) for all $j$.
\endroster
\enddefinition
\par
This is different from the standard definition of log-terminal
(resp. log-canonical) singularities,
we do not assume that $f^{-1}B$ is normal crossings.
\par
Note that if $B\geq 0$ then $f^*B\geq 0$.
Hence if $(Y,B)$ is log-terminal (resp. log-canonical) then
so is $(Y,0)$.
Moreover if $y\in\Supp(B)$ and $(Y,B)$ is log-canonical at $y$ then
$(Y,0)$ is log-terminal at $y$.
\par
\subhead 1.5\endsubhead
Here we recall some properties around $\delta_y$.
\par
\definition{Definition}
We define $\delta_y=-(Z-\Delta)^2$ and $\delta_{B,y}=-(Z-\Delta_B)^2$.
\enddefinition
\par
\proclaim{Proposition 3}{\rm ([KM, Theorem 1])}
\roster
\item $\delta_y=4$ if $y$ is smooth,
\item $\delta_y=2$ if $y$ is an RDP on $Y$,
\item $0<\delta_y<2$ if $(Y,0)$ is log-terminal at $y$,
\item $0\leq\delta_y\leq 2$ if $(Y,0)$ is log-canonical at $y$ but not smooth.
\endroster
\endproclaim
\demo{Proof}
(1) and (2) are clear.\par
For (3) and (4), $\delta_y\geq 0$ is obvious.
If $\delta_y=0$ then $Z=\Delta$, which implies log-canonical.
Hence if $(Y,y)$ is log-terminal at $y$ then $0<\delta_y$.
\par
On the other hand, $-(Z-\Delta)^2=-Z\cdot (Z+K_X)-K_X\cdot (Z-\Delta)$.
Since $f$ is minimal, $K_X\cdot (Z-\Delta)\geq 0$.
Furthermore $-Z\cdot (Z+K_X)=2-2p_a(Z)\leq 2$.
Hence we have $\delta_y\leq 2$.
If $\delta_y=2$, we have $p_a(Z)=2$ and
$K_X\cdot E_j=0$ for any $j$ with $e_j<z_j$.
If $(Y,0)$ is log-terminal at $y$, all $E_j$ satisfies $e_j<1\leq z_j$.
Hence $p_a(Z)=0$ and $Z\cdot K_X=0$ imply that $y$ is an RDP.
\qed
\enddemo
\par
\proclaim{Proposition 4}
Assume that $(Y,B)$ is log-canonical at $y$.
Then $\delta_{B,y}\leq \delta_y$.
Moreover the equality holds
if and only if $y\not\in\Supp(B)$.
\endproclaim
\par
\demo{Proof}
$$
\align
  \delta_y-\delta_{B,y}
  &= (Z-\Delta_B)^2-(Z-\Delta)^2 \\
  &= (\Delta-\Delta_B)\cdot (Z-\Delta_B+Z-\Delta) \\
  &= (f^{-1}B-f^*B)\cdot (Z-\Delta_B+Z-\Delta) \\
  &= f^{-1}B\cdot (Z-\Delta_B+Z-\Delta). \\
\endalign
$$
Since $Z-\Delta_B$ and $Z-\Delta$ is effective,
$\delta_{B,y}\leq\delta_y$.
\par
If $y\in\Supp(B)$ then $\Delta_B>\Delta$, therefore $(Y,0)$ is log-terminal.
Hence $f^{-1}B\cdot (Z-\Delta)>0$.
\qed
\enddemo
\par
\head 2. The main theorem\endhead
We assume the ground field is $\C$ throughout this paper.
We describe the statement of the main theorem first,
and then prove in section 3.
\par
\subhead 2.1\endsubhead
Let $Y$ be a projective normal surface over $\C$ and
$y$ be a fixed point on $Y$.
Let $M$ be a nef and big $\Q$-Weil divisor on $Y$ such that
$K_Y+\rup{M}$ is Cartier.
We set $B=\rup{M}-M$.
\par
Let $f\:(X,f^{-1}B)\to (Y,B)$ be the minimal resolution of the germ $(Y,y)$
if $y$ is singular, or the blowing up at $y$ if $y$ is smooth.
We define
$$
  \Delta_B=f^*(K_Y+B)-(K_X+f^{-1}B)=\sum e_jE_j,
$$
where $E_j$ are the exceptional curves lying over $y$.
\par
Let $Z=\sum z_jE_j$ be the fundamental cycle of $y$.
We define
$$
  \dmin
    =\min\{-(Z-\Delta_B+x)^2\mid
       \text{$x$ is an effective $f$-exceptional $\Q$-divisor}\}.
$$
\par
Since $-(Z-\Delta_B+x)^2$ is a quadric form defined by a negative definite
symmetric matrix of rational coefficients,
we have $\dmin$ is also rational and
there exists an effective $\Q$-divisor $x_0$ such that
$\dmin=-(Z-\Delta_B+x_0)^2$.
\par
We also define
$$
  \delta'=
    \cases
      1-\max\{e_1,e_n\}, &
        \text{if $(Y,B)$ is log-terminal at $y$ of type $A_n$,} \\
      & \text{where $E_1$ and $E_n$ are placed on the edge} \\
      & \text{of the chain of the dual graph} \\
      \text{any positive number}, &
        \text{if $(Y,B)$ is log-terminal at $y$ of type $D_n$} \\
      0, & otherwise
    \endcases
$$
\par
Let $\delta=\dmin$ if $(Y,B)$ is log-terminal at $y$
and $\delta=0$ otherwise.
We recall the main theorem introduced in section 0.
\par
\proclaim{Theorem 1}
If $M^2>\delta$ and $M\cdot C\geq\delta'$
for all curve $C$ passing through $y$
then $y\not\in\Bs{K_Y+\rup{M}}$.
\endproclaim
\par
Note that $\dmin\leq\delta_{B,y}\leq\delta_y\leq 4$.
Hence the bounds $\delta$ and $\delta'$ are effective.
\par
\subhead 2.2\endsubhead
We can get an easy corollary immediately, but we need some notation.
\par
\definition{Definition}
Assume $(Y,B)$ is log-terminal at $y$.
We define
$$
  \mu=\max\{t\geq 0\mid t(Z-\Delta)\leq f^*B\}.
$$
\enddefinition
\par
Since $(Y,0)$ is also log-terminal, $Z-\Delta>0$.
Hence $\mu$ is expressed as
$$
  \mu=\min\left\{\frac{b_j'}{z_j-a_j}\right\},
$$
where $\Delta=\sum a_jE_j$ and $f^*B-f^{-1}B=\sum b_j'E_j$.
Then we have $e_j=a_j+b_j'$ for all $j$.
\par
Note that $y\not\in\Supp(B)$ if and only if $\mu=0$.
If $y$ is smooth, $2\mu=\mult_y B$.
\par
\proclaim{Lemma 5}
If $(Y,B)$ is log-terminal at $y$ then $0\leq\mu<1$.
\endproclaim
\par
\demo{Proof}
If $\mu\geq 1$ then $b_j'\geq z_j-a_j$ for all $j$.
Hence $Z\leq\Delta_B$, that is contradiction.
\qed
\enddemo
\par
Let $x=(f^*B-f^{-1})-\mu (Z-\Delta)$ be an effective
$f$-exceptional divisor.
Since
$$
  (1-\mu)(Z-\Delta)=Z-\Delta-\sum b_j'E_j+x=Z-\Delta_B+x,
$$
we have
$$
  \dmin\leq (1-\mu)^2\delta_y.
$$
\par
We also have $\delta'\leq (1-\mu)\delta_y/2$.
Indeed, if $(Y,B)$ is not log-terminal of type $A_n$ at $y$, it is clear.
Hence we assume that $(Y,B)$ is log-terminal of type $A_n$ at $y$.
In this case we have $\delta_y=2-a_1-a_n$ by the following lemma.
Now we assume that $a_1\leq a_n$ by changing the indices.
Then we have
$$
  (1-\mu)\delta_y/2\geq (1-\mu)(1-a_n)=1-a_n-\mu(1-a_n)
  \geq 1-a_n-b_n'=1-e_n\geq\delta'.
$$
\par
\proclaim{Lemma 6}
$\delta_y=2-a_1-a_n$ if $(Y,B)$ is log-terminal of type $A_n$ at $y$.
In particular, if $n=1$ or $y$ is a smooth point then $\delta_y=2-2a_1$.
where the indices are taken in standard way.
\midinsert
\vskip 20pt
\includegraphics{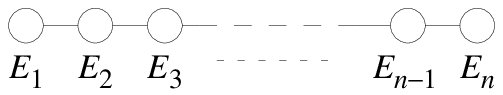}
\endinsert
\endproclaim
\par
\demo{Proof}
Since $y$ is rational, every $E_j$ is isomorphic to $\P^1$,
therefore $K_X\cdot E_j=-E_j^2-2$.
If $n\geq 2$ then $(Z-\Delta)\cdot E_j=0$ if $j\ne 1, n$ and
$(Z-\Delta)\cdot E_j=-1$ if $j=1$ or $n$.
Hence $-(Z-\Delta)^2=-(Z-\Delta)\cdot\sum (1-a_j)E_j=2-a_1-a_n$.
If $n=1$ then $(Z-\Delta)\cdot E_1=-2$.
Hence $-(Z-\Delta)^2=2-2a_1$.
\qed
\enddemo
\par
Now we have the following corollary of Theorem 1.
\par
\proclaim{Corollary 7}
If $M^2>(1-\mu)^2\delta$ and $M\cdot C\geq (1-\mu)\delta/2$
for all curves $C$ passing through $y$
then $y\not\in\Bs{K_Y+\rup{M}}$.
\endproclaim
\par
If $y$ is an RDP then $\Delta=0$ and $\delta_y=2$.
Thus Theorem 6 includes the result of Ein and Lazarsfeld ([EL Theorem 2.3]).
If $M$ is integral divisor, then $\mu=0$.
Hence this also includes the result of the author and Ma\c sek
(without the boundness of $C^2$, [KM Thoerm 2]).
\par
\head 3. A proof of the main theorem\endhead
In the case of that $(Y,B)$ is not log-terminal at $y$,
the proof is well known. (cf\. [ELM]).
So we assume that $(Y,B)$ is log-terminal at $y$.
\par
\subhead 3.1\endsubhead
Before starting the proof of Theorem 1,
we reduce the problem to the situation where $(Y,y)$ has no other singularity.
\par
\proclaim{Lemma 8}{\rm ([Ma\c s, Lemma 10])}
We may assume that $Y-\{y\}$ is smooth.
\endproclaim
\par
\demo{Proof}
If $Y$ has other singularities,
then take $g\:S\to Y$ be a simultaneous resolution of all the singularities
of $Y$ except $y$.
\par
Let $M'=g^*M$ and $y'=g^{-1}(y)$ be a point on $S$.
Since $M$ is nef and big, $M'$ is also nef and big.
Since $S-\{y'\}$ is smooth and $\rup{M'}$ is integral coefficients,
$K_S+\rup{M'}$ is Cartier except $y'$.
Since $g$ is isomorphism on a neighbourhood of $y'$,
we have $K_S+\rup{M'}$ is also Cartier divisor on $S$.
\par
Now we show that if the theorem is true for $S, M', y'$
then it is true for $Y,M,y$.
\par
Since $M'=g^*M$ we have $(M')^2>\delta$ and $M'\cdot C'\geq\delta'$
for all curves $C'$ passing through $y'$.
If the theorems are true for $S, M', y'$ then there exists a section
$s'\in H^0(S,K_S+\rup{M'})$ such that $s'(y')\ne 0$.
\par
Let $\Delta'=g^*K_Y-K_S$ be an effective divisor on $S$.
Then we have
$$
  K_S+\rup{M'} = \rup{K_S+g^*M} =
  g^*(K_Y+\rup{M})-\rdown{\Delta'+g^*B}.
$$
Since all coefficients of $f^{-1}B$ are less than 1,
$N=\rdown{\Delta'+g^*B}$ is $g$-exceptional divisor and $y'\not\in\Supp(N)$.
Multiplying $s'$ by the global section on $\O_S(N)$,
we find a section $s\in H^0(S, g^*(K_Y+\rup{M}))$
such that $s(y')\ne 0$.
Since $Y$ is normal, $s$ corresponds to a global section of
$\O_Y(K_Y+\rup{M})$ which does not vanish at $y$.
\qed
\enddemo
\par
\subhead 3.2\endsubhead
Now we assume that $Y$ is smooth except $y$.
We assume $(Y,B)$ is log-terminal at $y$ throughout this paper.
First we introduce the following lemma.
\par
Let $f\:X\to (Y,y)$ be a birational morphism from smooth surface
to a germ of a normal surface singularity.
Let $\Gamma$ be an $f$-exceptional $\Q$-Weil divisor on $X$ and
let $\gamma=-\Gamma^2>0$.
Let $M$ be a nef and big $\Q$-Weil divisor on $Y$.
\proclaim{Lemma 9}
Assume $M^2>\gamma$.
Then there exists an effective $\Q$-Weil divisor $D$ on $Y$ such that
\roster
\item $D\equiv M$,
\item $f^*D>\Gamma$,
\endroster
\endproclaim
\par
\demo{Proof}
Since $M^2>\gamma$,
we may assume $M^2>(1+\sigma)^2\gamma=-((1+\sigma)\Gamma)^2$
for very small rational number $0<\sigma\ll 1$.
Thus we replace $\Gamma'=(1+\sigma)\Gamma$ by $\Gamma$,
it is enough to show that (1) and (2${}'$) $f^*D\geq\Gamma$.
\par
Since $(f^*M-\Gamma)^2>0$ and $f^*M\cdot (f^*M-\Gamma)>0$,
the $\Q$-divisor $f^*M-\Gamma$ is in the positive cone of $X$.
Hence $f^*M-\Gamma$ is big and there is a member
$T\in |n(f^*M-\Gamma)|$ for sufficiently large and divisible $n$.
Then we set $T'=\frac1n T+\Gamma\equiv f^*M$.
\par
Let $D=f_*T'$.
\roster
\item $D$ is effective because $T'$ is effective.
\item Let $G=f^*D-T'$ be an $f$-exceptional $\Q$-divisor.
      Since $G^2=(f^*D-T')\cdot G=(f^*D-f^*M)\cdot G=0$, we have $f^*D=T'$.
\item Since $f^*D=T'\equiv f^*M$, we have $D\equiv M$.
\item $f^*D=T'\geq \Gamma$.
\endroster
\qed
\enddemo
\par
\subhead 3.3\endsubhead
We come back to the proof of Theorem 1.
\par
Since $M^2>\dmin$,
there exists an effective $\Q$-Weil divisor $D$
such that
\roster
\item $D\equiv M,$
\item $f^*D>Z-\Delta_B+x$ where $x=\sum x_iE_i\geq 0$ gives the minimum
      $\dmin$,
\endroster
by Lemma 11.
\par
Let $D=\sum d_iC_i$, $B=\sum b_iC_i$ and $D_i=f^{-1}C_i$.
Let $f^*D=\sum d_iD_i+\sum d_j'E_j$ and $f^*B=\sum b_iD_i+\sum b_j'E_j$.
Note that we take the sum $\sum b_iD_i$ and $\sum d_iD_i$ commonly,
so some of $b_i$ and $d_i$ may be zero.
\par
We define a rational number $c$ as follows.
$$
  c =\min\left\{\frac{1-e_j}{d_j'},\frac{1-b_i}{d_i}\Bigm|
    \vcenter{\hsize=2.8in\parindent=0pt
      for all $j$ such that $f(E_j)=\{y\}$ and \hfil\break
      for all $i$ such that $D_i$ meets $f^{-1}(y)$ and $d_i>0$}
    \right\}.
$$
Since $(Y,B)$ is log-terminal at $y$ and $b_i<1$ for all $i$,
we have $c>0$.
Since $f^*D>Z-\Delta_B+x$, we have 
$d_j'>z_j-e_j+x_j$.
Hence we have
$$
  d_j'+e_j> z_j+x_j\geq 1.
$$
Thus we have $c<1$.
\par
Let $R=f^*M-cf^*D\equiv (1-c)f^*M$.
Hence $R$ is nef and big.
Note that
$$
\align
  \rup{R}
    &= \rup{K_X+R}-K_X=\rup{K_X+f^*M-cf^*D}-K_X \\
    &= \rup{f^*(K_Y+\rup{M})-f^*B-cf^*D-\Delta}-K_X \\
    &= f^*(K_Y+\rup{M})-K_X-\rdown{cf^*D+f^*B+\Delta} \\
    &= f^*M+f^*B+\Delta-\rdown{cf^*D+f^*B+\Delta} \\
    &= R+\{cf^*D+f^*B+\Delta\}
\endalign
$$
Hence we have
$$
\align
  K_X+\rup{R}
    &= f^*(K_Y+\rup{M})-\rdown{cf^*D+f^*B+\Delta} \\
    &= f^*(K_Y+\rup{M})-\sum\rdown{cd_i+b_i}D_i-\sum\rdown{cd_j'+e_j}E_j.
\endalign
$$
By the definition of $c$, we have $0<cd_j'+e_j\leq 1$ for all $j$
and $0<cd_i+b_i\leq 1$ for any $i$ with $\Supp D_i\cap f^{-1}(y)\ne\emptyset$.
Let $\sum\rdown{cd_i+b_i}D_i=A+N$ where $\Supp(N)\cap f^{-1}(y)=\emptyset$
and $A=D_1+\cdots +D_t$ where $D_1,\dots ,D_t$ meets $f^{-1}(y)$.
Let $E=\sum\rdown{cd_j'+e_j}E_j$.
\par
By the minimality of $c$, at least one of $E$ and $A$ is non-zero.
Each component $E_j$ in $E$ and each component $D_i$ in $A$ have
coefficients $1$ in $R$,
they do not appear in the fractional part of $R$.
Since $D_i\in A$ meets $f^{-1}(y)$,
we have $\rup{R}\cdot D_i\geq R\cdot D_i=(1-c)f^*M\cdot D_i\geq (1-c)\delta'$.
\par
\proclaim{Lemma 10}
If $A\ne 0$ then $(Y,f_*A)$ is log-canonical at $y$.
Moreover,
let $\Gamma$ be the dual graph of the union of the exceptional locus $E$
and the strict transforms $A$ of the $f_*A$.
Then one of the following holds, where \kern10pt\lower-3pt\hbox{\maru j} and
\kern10pt\lower-3pt\hbox{\maru z} indicate
prime components of $E$ and $f_*A$, respectively.
\roster
\item $f_*A$ is irreducible and non-singular, $(Y,y)$ is a cyclic quotient
      singularity, and $\Gamma$ is as follows:
\vskip 20pt
\includegraphics{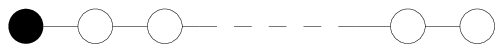}
\item $f_*A$ is irreducible and non-singular, 
      $(Y,y)$ is a quotient singularity and $\Gamma$ is as follows:
\vskip 40pt
\includegraphics{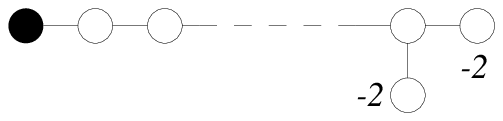}
\item $f_*A$ has two prime components which are non-singular and
      intersect transversally, $(Y,y)$ is a cyclic quotient singularity
      and $\Gamma$ is as follows:
\vskip 20pt
\includegraphics{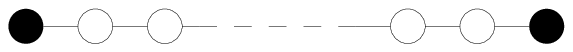}
\endroster
\endproclaim
\par
\demo{Proof}
Since
$$
  f^*(K_Y+f_*A)-K_X-A\leq f^*(cD+B)-A+\Delta\leq \sum (cd_j'+e_j)E_j,
$$
$(Y,f_*A)$ is log-canonical at $y$.
These are classified as in [Al] and [Ky], they are only above 3 cases.
\qed
\enddemo
\par
We come back to the proof of Theorem 1.
We consider two cases according to whether $E$ is zero or not.
\par
\medskip
Case 1: $E\ne 0$.
\par
We recall that the following Kawamata-Viehweg vanishing theorem.
\par
\proclaim{Claim 11 {\rm (cf\. [Sa, Lemma 5])}}
Let $X$ be a smooth projective surface over $\Bbb C$ and
let $R$ be a nef and big $\Q$-divisor on $X$.
Let $D_1,\dots D_t$ be distinct irreducible curves such that
they does not appear in the fractional part of $R$ and
$R\cdot D_i>0$ for any $i$.
Then
$$
  H^1(X,K_X+\rup{R}+D_1+\cdots +D_t)=0.
$$
\endproclaim
\par
If $t>0$ then $A\ne 0$ and $y$ is of type $A_n$ or $D_n$ by the above lemma.
Hence $R\cdot D_i=(1-c)f^*M\cdot D_i>0$ by assumption.
So it applies
$$
  H^1(X, K_X+\rup{R}+A)=H^1(X, K_X+\rup{R}+D_1+\cdots +D_t)=0.
$$
That is, $H^1(X,f^*(K_Y+\rup{M})-E-N)=0$.
Hence
$$
  H^0(X, f^*(K_Y+\rup{M})-N)\to H^0(E, (f^*(K_Y+\rup{M})-N)|_E)
$$
is surjective.
Then we get a global section $s\in H^0(X,f^*(K_Y+\rup{M})-N)$
which does not vanish anywhere on $f^{-1}(y)$.
Since $\Supp(N)\cap f^{-1}(y)=\emptyset$, multiplying $s$ by
global section of $\O(N)$, we find a global section $t\in H^0(Y,K_Y+\rup{M})$
such that $t(y)\ne 0$.
\par
\medskip
Case 2: $E=0$.
\par
In this case, $A\ne 0$.
\par
Note that
$f^*(K_Y+f_*A)-K_X-A\leq\sum (cd_j'+e_j)E_j$.
Since $E=0$, we have $cd_j'+e_j <1$ for all $j$.
Hence $(Y,f_*A)$ is log-terminal with reduced boundary $f_*A$,
so it is of type $A_n$ by Lemma 12.
Moreover $t=1$, namely $A=D_1$.
\par
Note that if $y$ is a singular point of type either $D_n$ or $E_n$, 
This case 2 never occurs.
\par
By Kawamata-Viehweg vanishing theorem we have
$$
  H^1(X, K_X+\rup{R})=0.
$$
That is, $H^1(X,f^*(K_Y+\rup{M})-D_1-N)=0$.
Hence
$$
  H^0(X, f^*(K_Y+\rup{M})-N)\to H^0(D_1, (f^*(K_Y+\rup{M})-N)|_{D_1})
$$
is surjective.
Hence it is enough to find a global section
$s\in H^0(D_1, (f^*(K_Y+\rup{M})-N)|_{D_1})$ which does not vanish at
$p\in \Supp (D_1)\cap f^{-1}(y)$.
Indeed, if there is such $s$, there is also a global section
$s'\in H^0(X,f^*(K_Y+\rup{M})-N)$ such that $s'(p)\ne 0$.
By multiplying a section of $\O(N)$, we find a desired section
$t\in H^0(Y,K_Y+\rup{M})$ which does not vanish at $y$.
\par
Since
$$
  (f^*(K_Y+\rup{M})-N)|_{D_1}=(K_X+\rup{R}+D_1)|_{D_1}
  =K_{D_1}+\rup{R}|_{D_1},
$$
$\rup{R}\cdot D_1>1$ implies
$\rup{R}\cdot D_1\geq 2$ because it is integer.
This yield a section
$s\in H^0(D_1, K_{D_1}+\rup{R}|_{D_1})$
which does not vanish at $p$ by the theorem in [H].
Therefore it is enough to prove $\rup{R}\cdot D_1>1$.
\par
Since
$$
  \rup{R}=R+\sum_{i\geq 2}\{cd_i+b_i\}D_i+\sum\{cd_j'+e_j\}E_j
$$
and $E=\sum \rdown{cd_j'+e_j}E_j=0$,
we have
$$
  \rup{R}\cdot D_1\geq R\cdot D_1+\sum (cd_j'+e_j)E_j\cdot D_1.
$$
\par
Since $y\in f_*D_1$,
we have $R\cdot D_1=(1-c)f^*M\cdot D_1\geq (1-c)\delta'$
and $E_j\cdot D_1=1$ for $j=1$ or $j=n$.
By changing the index of $\{E_j\}$, we may assume that $e_1\leq e_n$.
\par
Case 2-1: $D_1$ meets $E_n$.
\par
Recall that $\delta'=1-\max\{e_1,e_n\}=1-e_n$
since we take $e_1\leq e_n$
In this case,
$$
\align
  \rup{R}\cdot D_1
  &\geq (1-c)\delta'+(cd_n'+e_n) \\
  &= (1-c)(1-e_n)+cd_n'+e_n \\
  &= 1+c(d_n'+e_n-1).
\endalign
$$
Since $d_j'>z_j-e_j+x_j$ for all $j$,
we have $\rup{R}\cdot D_1>1$
because $y$ is of type $A_n$.
\par
\medskip
Case 2-2: $D_1$ meets $E_1$.
\par
Let $A=A(w_1,\dots,w_n)=(-E_i\cdot E_j)_{i,j}$ be the (positive definite)
intersection matrix of type $A_n$ where $w_j=-E_j^2$.
Let $a(w_1,\dots,w_n)=\det A(w_1,\dots,w_n)$.
Let $\bold b={}^t(b_1',\dots,b_n')$ where $f^*B=f^{-1}B+\sum b_j'E_j$.
Let $L_i$ be an irreducible reduced curve on $Y$ such that
$f^{-1}L_i\cdot E_i=1$ for an $i$ and $f^{-1}L_i\cdot E_j=0$ for $j\ne i$.
Let $f^*L_i=f^{-1}L_i+\sum c_{ij}E_j$.
\par
Since $b_j'$ are defined by linear equation
$A\bold b={}^t(-f^{-1}B\cdot E_1,\dots,-f^{-1}B\cdot E_n)$,
we have $b_j'=\sum (f^{-1}B\cdot E_i)c_{ij}$.
Hence it is enough to calculate $c_{ij}$ for $b_j'$.
\par
We define $a()=1$ for convenience.
\par
\proclaim{Proposition 12}
\roster
\item
$$
  1-a_i=\frac{a(w_1,\dots,w_{i-1})+a(w_{i+1},\dots,w_n)}{a(w_1,\dots,w_n)},
$$
\item
$$
  c_{ij}=\frac{a(w_1,\dots,w_{i-1})a(w_{j+1},\dots,w_n)}{a(w_1,\dots,w_n)},
         \text{ if }i\leq j.
$$
Moreover we have $c_{ij}=c_{ji}$.
\endroster
\endproclaim
\par
\demo{Proof}
Let $\bold a={}^t(a_1,\dots,a_n)$,
$\bold w={}^t(-w_1+2,\dots,-w_n+2)$,
$\bold c={}^t(c_{i1},\dots,c_{in})$ and
$\bold e_i={}^t(0,\dots,0,1,0,\dots,0)$ where the only $i$-th element is 1.
\par
Let $\widetilde{A}_{ij}$ be the $(i,j)$-component of $A^{-1}$.
If $i<j$ then we have
$$
  \widetilde{A}_{ij}=
    \frac1{\det A}(-1)^{i+j}a(w_1,\dots,w_{i-1}) a(w_{j+1},\dots,w_n).
$$
Also we have
$$
  \widetilde{A}_{ii}=
    \frac1{\det A}a(w_1,\dots,w_{i-1}) a(w_{i+1},\dots,w_n).
$$
By the wipe-out method, we have
$$
  a(w_1,\dots,w_j)=w_j a(w_1,\dots,w_{j-1})-a(w_1,\dots,w_{j-2}).
$$
\par
(1) Since $\bold a=-A^{-1}\bold w$, calculate the right hand side
using above equalities.
\par
(2) Since $\bold c=-A^{-1}\bold e_j$, calculate the right hand side
using above equalities.
\par
(For more detail, see [Kt]).\qed
\enddemo
\par
Especially we have
$$
\align
  1-a_1 &= (1+a(w_2,\dots,w_n))/a(w_1,\dots,w_n), \\
  1-a_n &= (1+a(w_1,\dots,w_{n-1}))/a(w_1,\dots,w_n), \\
  c_{1n} &= c_{n1}=1/a(w_1,\dots,w_n), \\
  c_{1j} &= a(w_{j+1},\dots,w_n)/a(w_1,\dots,w_n).
\endalign
$$
Suppose $n\geq 2$.
Since $a(w_j,\dots,w_n)=w_j a(w_{j+1},\dots,w_n)-a(w_{j+2},\dots,w_n)$
and $w_j\geq 2$ for all $j$,
we have
$$
\align
  a(w_j &,\dots,w_n)-a(w_{j+1},\dots,w_n) \\
  &= (w_j-1)a(w_{j+1},\dots,w_n)-a(w_{j+2},\dots,w_n) \\
  &\geq a(w_{j+1},\dots,w_n)-a(w_{j+2},\dots,w_n)
   \geq a(w_n)-1 = w_n-1>0.
\endalign
$$
Hence we have $c_{11}>c_{12}>\dots>c_{1n}$ and samely we have
$c_{n1}<c_{n2}<\dots<c_{nn}$.
\par
\proclaim{Lemma 13}
Let $P$ be an effective $\Q$-divisor on $Y$.
Let $f^*P=f^{-1}P+\sum p_jE_j$.
Then we have $p_n\leq a(w_1,\dots,w_{n-1}) p_1$.
\endproclaim
\par
\demo{Proof}
If we set $q_i=f^{-1}P\cdot E_i$, we have $p_j=\sum q_ic_{ij}$ for all $j$.
Hence
$$
  p_1=\sum q_ic_{i1}\geq \sum q_ic_{n1}
     =\frac1{a(w_1,\dots,w_n)}\sum q_i,
$$
and
$$
  p_n=\sum q_ic_{in}\leq \sum q_ic_{nn}
     =\frac{a(w_1,\dots,w_{n-1})}{a(w_1,\dots,w_n)}\sum q_i.
$$
Therefore $p_n\leq a(w_1,\dots,w_{n-1}) p_1$.
\qed
\enddemo
\par
We come back to the proof of Theorem 1.
\par
Let $f^*C_1=D_1+\sum_{j=1} c_j E_j$.
Note that $c_j=a(w_{j+1},\dots,w_n)/\alpha$,
where $\alpha=a(w_1,\dots,w_n)$.
Let $y_{D,j}=d_j'-d_1 c_j$, $y_{B,j}=b_j'-b_1 c_j$ and $y_j=cy_{D,j}+y_{B,j}$.
Note that $a_1+c_1=1-1/\alpha$ by Proposition 10.
\par
Since $cd_1+b_1=1$, we have
$$
  cd_1'+b_1'
  = c(y_{D,1}+d_1c_1)+y_{B,1}+b_1c_1
  = c_1+cy_{D,1}+y_{B,1}
$$
\par
In this case,
$$
\align
  \rup{R}\cdot D_1
  &\geq (1-c)\delta'+cd_1'+e_1 \\
  &= (1-c)(1-e_n)+cd_1'+e_1 \\
  &= (1-c)(1-e_n)+c_1+a_1+(cy_{D,1}+y_{B,1})\\
  &= (1-c)(1-e_n)+1-\frac1{\alpha}+y_1.
\endalign
$$
Since $E=0$, we have $cd_1'+e_1<1$.
Hence $cd_1+b_1'<1-e_1+b_1'=1-a_1$.
Thus $y_1=c(d_1'-d_1c_1)+(b_1'-b_1c_1)<1-a_1-c_1=1/\alpha$
by Proposition 10.
\par
\proclaim{Claim 14}
$$
  (1-c)(1-e_n)>\frac{a(w_1,\dots,w_{n-1})}{\alpha}-y_n.
$$
\endproclaim
\par
\demo{Proof}
By the choice of $D$, we have $d_n'>1-e_n+x_n\geq 1-a_n-b_n'$.
Hence
$$
  (d_n'-1+a_n+b_n')\frac{c}{1-a_n}
  >0=\frac{cd_1+b_1-1}{1+a(w_1,\dots,w_{n-1})},
$$
since $cd_1+b_1=1$.
We set $\alpha'=a(w_1,\dots,w_{n-1})$ for convenience.
Then we have
$$
  \left(
    (d_n'-1+a_n+b_n')\frac1{1-a_n}-\frac{d_1}{1+\alpha'}
  \right)c > \frac{b_1-1}{1+\alpha'}.
$$
Since $(1-a_n)\alpha=1+\alpha'$ and $d_n'=d_1c_n+y_{D,n}=d_1/\alpha+y_{D,n}$,
the left-hand-side is equal to
$$
  \left(\frac{d_n'}{1-a_n}-1+\frac{b_n'}{1-a_n}-\frac{d_1}{1+\alpha'}\right)c
  =\left(\frac{y_{D,n}}{1-a_n}+\frac{b_n'}{1-a_n}-1\right)c.
$$
\par
On the other hand, the right-hand-side is equal to
$$
\align
  \frac{b_1-1}{1+\alpha'}
  &= \frac{b_1+\alpha y_{B,n}}{1+\alpha'}-\frac{1+\alpha y_{B,n}}{1+\alpha'}\\
  &= \frac{b_n'}{1-a_n}-\frac{1+\alpha y_{B,n}}{1+\alpha'} \\
  &= \frac{b_n'}{1-a_n}-1+\frac{\alpha'-\alpha y_{B,n}}{1+\alpha'}. \\
\endalign
$$
\par
Hence we have the inequality
$$
  \left(\frac{y_{D,n}}{1-a_n}-1+\frac{b_n'}{1-a_n}\right)c
  >\frac{b_n'}{1-a_n}-1+\frac{\alpha'/\alpha-y_{B,n}}{1-a_n}.
$$
Thus we have
$$
\align
  (1-c)\left(1-\frac{b_n'}{1-a_n}\right) 
    &> \frac{\alpha'/\alpha-y_{B,n}-cy_{D,n}}{1-a_n} \\
  (1-c)(1-a_n-b_n') &> \frac{\alpha'}{\alpha}-y_n.
\endalign
$$
\qed
\enddemo
\par
By this claim, we have
$$
  \rup{R}\cdot D_1>1+\frac{\alpha'-1}{\alpha}+y_1-y_n.
$$
Since $f^*(c(D-C_1)+(B-C_1))=\sum y_jE_j$,
we have $y_n\leq\alpha' y_1$ by Lemma 11.
Hence we have
$$
  \frac{\alpha'-1}{\alpha}+y_1-y_n
  \geq\left(\alpha'-1)(\frac1{\alpha}-y_1\right).
$$
Since $\alpha'=a(w_1,\dots,w_{n-1})\geq 1$ and $y_1<1/\alpha$,
we have $\rup{R}\cdot D_1>1$.
\qed
\par

\enddocument